
\documentstyle{amsppt}
\magnification=\magstep1
\parindent=1em
\baselineskip 15pt
\hsize=12.3 cm
\vsize=18.5 cm
\NoRunningHeads
\NoBlackBoxes
\pageno=1


\def\supp{\text{supp }}
\def\tr{\text{tr }}

\def\dim{\text{dim }}



\topmatter

\title Uniqueness of unconditional bases in $c_0-$products
\endtitle
\author
P.G. Casazza
and N.J. Kalton
\endauthor
\address
Department of Mathematics,
University of Missouri,
Columbia, Mo.  65211, U.S.A.
\endaddress
\email pete\@math.missouri.edu, nigel\@math.missouri.edu
\endemail
\thanks The first author was  supported  by  NSF Grant 9706108; the
second author was supported by  Grant DMS-9500125 \endthanks
\abstract
We give counterexamples to a conjecture of Bourgain, Casazza,
Lindenstrauss and Tzafriri that if $X$ has a unique unconditional basis
(up to permutation) then $c_0(X)$ also has a unique unconditional basis.
We also give some positive results including a simpler proof that
$c_0(\ell_1)$ has a unique unconditional basis and a proof that
$c_0(\ell_{p_n}^{N_n})$ has a unique unconditional basis when
$p_n\downarrow 1$, $N_{n+1}\ge 2N_n$ and $(p_n-p_{n+1})\log N_n$ remains
bounded.

\endabstract

\endtopmatter

\heading{1. Introduction}\endheading

A Banach space $X$ is said to have a {\it unique unconditional basis}
(or more precisely, a unique unconditional basis up to permutation) if
it has an unconditional basis and  if whenever $(u_n)$ and $(v_n)$ are
two normalized unconditional bases of $X$, then there is permutation
$\pi$ of
$\Bbb N$ such that $(v_n)$ and $(u_{\pi(n)})$ are equivalent.  Since
unconditional bases correspond to discrete or atomic order-continuous
lattice structures on $X$, this can be reworded as a statement that such
a lattice-structure is essentially unique.

The earliest examples of Banach spaces with unique unconditional
bases are
$c_0,\ell_1$
(\cite{10}) and $\ell_2$ (\cite{9}). It was shown by Lindenstrauss and
Zippin \cite{12} that amongst spaces with symmetric bases this is the
complete list. Later Edelstein and Wojtaszczyk
showed that direct sums of these spaces also have unique unconditional
bases. All these results can be found in \cite{11}. In \cite{3} the
authors attempted a complete classification
and showed that the spaces $c_0(\ell_1),c_0(\ell_2),\ell_1(c_0)$ and
$\ell_1(\ell_2)$ all have unique unconditional bases while
$\ell_2(\ell_1)$ does not.  They also found an unexpected additional
space, 2-convexified Tsirelson (see \cite{5} for the definition) with a
unique unconditional basis.
Recently, the authors found a new approach to this type of problem and
were  able to add some more spaces, including Tsirelson space (see
\cite{5})
itself and certain Nakano spaces \cite{4} (as pointed out in \cite{4}
some spaces considered by Gowers \cite{8} provide further examples)
; we also showed
that uniqueness of the unconditional basis need not be inherited  by a
complemented subspace.

This note is motivated by a question raised in \cite{3}.  They asked
if $c_0(X)$ has a unique unconditional basis whenever $X$ does.  The
idea here is that if this and the corresponding dual result for
$\ell_1-$products holds then one could iterate the results in
\cite{3} to produce examples such as $c_0(\ell_1(c_0(\ell_1)))$ and
so on.

Unfortunately as we show below in Section 4, the answer to this question
is negative and Tsirelson space $T$ or its 2-convexified version both
produce counterexamples.  However, we show how our approach in \cite{4}
can be used for $c_0-$products.  We give a much shorter
proof
(Theorem 3.3) of the fact that $c_0(\ell_1)$ has a unique unconditional
basis; the original proof of this result in \cite{3} is extremely
technical. We show by the same techniques
(Theorem
3.4) that examples of the type
$c_0(\ell_{p_n}^{N_n})$ where $p_n\downarrow 1$, $N_{n+1}\ge 2N_n$ and
$(p_n-p_{n+1})\log N_n$ remains bounded must also have unique
unconditional bases.

In Section 4, we also use the same techniques to show that for certain
right-dominant spaces $X$, as introduced in \cite{4}, such as Tsirelson
space
$T$, any unconditional basis of $c_0(X)$ must be equivalent to a
{\it subset} of the canonical basis (Theorem 4.1).  Nevertheless we show
that the unconditional basis of $c_0(T)$ is {\it not} unique as already
remarked.

We conclude this section with a few remarks on terminology and
assumptions.  We will frequently index unconditional bases and basic
sequences by an unordered countable index set $\Cal N$ which need not be
the natural numbers $\Bbb N.$  We will assume that any unconditional
basic
sequence $(u_n)_{n\in\Cal N}$ is semi-normalized, i.e. $0<\inf_{n\in\Cal
N}\|u_n\|\le\sup_{n\in\Cal N}\|u_n\|<\infty.$  We will say that two
unconditional basic sequences $(u_n)_{n\in\Cal N}$ and $(v_n)_{n\in\Cal
N'}$ are {\it equivalent} if there there is a bijection $\pi:\Cal
N\to\Cal N'$ so that $(u_n)_{n\in\Cal N}$ and $(v_{\pi(n)})_{n\in\Cal N}$
are equivalent.

An unconditional basic sequence $(u_n)_{n\in\Cal N}$ in $X$ is {\it
complemented} if there is a bounded projection $P:X\to [u_n]_{n\in\Cal
N}.$   If $(u_n)_{n\in\Cal N}$ is an unconditional basis of $X$ and
$(v_n)_{n\in\Cal N'}$ is an unconditional basic sequence of the form
$v_n=\sum_{k\in\Cal A_n}a_ku_k$ where the sets $(\Cal A_n)_{n\in\Cal N'}$
are disjoint subsets of $\Cal N$ we say that $(v_n)_{n\in\Cal N'}$ is
{\it disjoint} with respect to $(u_n)_{n\in\Cal N}.$  If $(v_n)$ is a
complemented disjoint sequence then it may be shown that there is a
projection $Px=\sum_{n\in\Cal N'}v_n^*(x)v_n$ where each $v_n^*\in X^*$
is of the form $v_n^*=\sum_{k\in\Cal A_n}b_ku_k^*$ where
$(u_k^*)_{k\in\Cal N}$ are the sequence of biorthogonal functions for
$(u_k)_{k\in\Cal N}.$

It will be convenient to represent a space $X$ with unconditional basis
$(u_n)_{n\in\Cal N}$ as a sequence space modelled on the index set $\Cal
N$, identifying $\sum_{k\in\Cal N}a_ku_k$ with the function $f:\Cal N\to
\Bbb R$ given by $f(k)=a_k.$  This identifies $X$ as a discrete Banach
lattice and allows us to use functional notation.  The canonical basis of
a sequence space $X$ modelled on $\Cal N$ is denoted by $(e_n)_{n\in\Cal
N}.$

If $(u_n)_{n\in\Cal N}$ is an unconditional basis for $X$ and $N$ is a
natural number we denote by $(u_n)_{n\in\Cal N}^N$ the naturally induced
unconditional basis of $X^N$ (the direct sum of $N$ copies of $X$).

For future reference we note here that our techniques depend critically
on the following result, proved in Theorem 3.5 of \cite{4}:

\proclaim{Theorem 1.1}Suppose $X$ is a Banach space with an unconditional
basis $(u_n)_{n\in\Cal N}$ which does not contain uniformly complemented
copies of
$\ell_2^n$
(i.e. is not sufficiently Euclidean). Suppose $(v_n)_{n\in\Cal B}$ is a
complemented unconditional basic sequence in $X$.  Then there is a
integer $N$ and a complemented disjoint sequence $(w_n)_{n\in\Cal B}$ in
the basis
$(u_n)_{n\in\Cal N}^N$ such that $(v_n)_{n\in\Cal B}$ is equivalent to
$(w_n)_{n\in\Cal B}.$\endproclaim

\heading{2. A criterion for an $\ell_1-$ or $c_0-$product to be
sufficiently Euclidean}\endheading

The aim of this section is to establish criteria for $c_0$-product
to contain  uniformly complemented copies of $\ell_2^n$ so that
we can apply Theorem 1.1.

If $X$ is a Banach space we will say that $X$ has property $P(k,M)$ where
$k\in \Bbb N$ and $M\ge 1$ if whenever $S:\ell_2^k\to X$ and
$T:X\to\ell_2^k$ are operators satisfying $TS=I_{\ell_2^k}$ then
$\|S\|\|T\|\ge M.$  We will say that a sequence of Banach spaces
$(X_j)_{j=1}^{\infty}$ has property $P(k,M)$ if each $X_j$ has property
$P(k,M).$  A Banach space $X$ (respectively a sequence of Banach spaces
$(X_j)_{j=1}^{\infty}$) is {\it sufficiently Euclidean} if there exists
$M$ so
that $X$ (respectively $(X_j)_{j=1}^{\infty})$ fails $P(k,M)$ for every
$k\in\Bbb N.$

We recall that if $H$ is a finite-dimensional Hilbert space and $A:H\to
X$ is any linear map then the $\ell-$norm of $A$ is given by:
$$ \ell(A) =\bold E(\|\sum_{i=1}^mg_iAe_i\|^2)^{1/2}$$ where
$(e_,\ldots,e_m)$ is any orthonormal basis of $H$ and $(g_1,\ldots,g_m)$
is a sequence of independent normalized Gaussian random variables.  See
\cite{15}. If $S$ is an operator on a Banach space $X$ and $E$ is a
closed subspace of $X$ we denote by $S_E$ the restriction of $S$ to $E.$

\proclaim{Lemma 2.1}There exists a universal constant $c>0$
with the following property.  Suppose
$H$ is an
$n-$dimensional Hilbert space and
$X$
is any Banach space.  Suppose $S:H\to X$ is any   operator with
$\|S\|\le
1$.  Then there
is a subspace $E$ of $H$ with $\dim E\ge c\ell(S)^2$ so that $\|S_E\|\le
3\ell(S)n^{-1/2}.$
\endproclaim

\demo{Proof}It will suffice to prove this for $S$ one-to-one, since
the result then follows by a simple perturbation argument. Let
$\mu$ be normalized invariant measure on the surface
of the sphere in $\ell_2^n.$  Consider the norm $\xi\to \|S\xi\|$; this
satisfies $\|S\xi\|\le \|\xi\|$ for all $\xi.$  We
use
Theorem 4.2 of \cite{13} (p.12).  If $M_r$ is a median value of the
norm
$\|S\xi\|$ then
$$ M_r \le \sqrt 2 (\int \|S\xi\|^2d\mu)^{1/2}=\sqrt{2/n}\ell(S).$$
\enddemo

\proclaim{Lemma 2.2}Suppose $X$ is a Banach space with property
$P(k,M).$
Suppose $H$ is an $n-$dimensional Hilbert space and $S:H\to X$ and $T:
X\to H$ are bounded operators with $\|T\|\le 1.$  Then
$$ |\tr(TS)| \le Cn^{1/2}\max(M^{-1}\ell(S),k^{1/2}\|S\|)$$
for some universal constant $C.$

\endproclaim

\demo{Proof}Suppose $1\le j\le n$ and that $s_j$ is the $j$th. singular
value of $TS.$  We can restrict to a subspace $H_j$ of dimension $j$ so
that $\|TS\xi\|\ge s_j\|\xi\| $ for all $\xi\in H_j.$

Assume $s_j>0.$  Then by the preceding Lemma 2.1 there is a subspace $E$
of $H_j$ so that $\dim E\ge c\ell(S_{H_j})^2$ and $\|S_E\|\le
3\ell(S_E)j^{-1/2}.$

Suppose $\dim E\ge k.$  Then we must have $\|T\|\|S_E\|\ge Ms_j$ so that
$\ell(S_E) \ge \frac13 Ms_jj^{1/2}.$  This implies that $s_j \le
3M\ell(S)j^{-1/2}.$  If $\dim E<k$ then $$\ell(S_{H_j})^2\|S_{H_j}\|^{-2}
\le c^{-1}k$$
and so $\ell(TS_{H_j}) \le c^{-1/2}k^{1/2}\|S\|.$  From this we deduce
that
$j^{1/2}s_j \le c^{-1/2}k^{1/2}\|S\|$ or $s_j\le
c^{-1/2}k^{1/2}j^{-1/2}\|S\|.$ Combining we obtain that
$$ s_j \le \max(3M^{-1}\ell(S),c^{-1/2}k^{1/2}\|S\|)j^{-1/2}.$$

Now $$|\tr (TS)|\le \sum_{j=1}^ns_j \le
Cn^{1/2}\max(M^{-1}\ell(S),k^{1/2}\|S\|)$$
for some universal constant $C.$ \enddemo

\proclaim{Lemma 2.3}There is a universal constant $C$ so that if $X$ is a
Banach space with property $P(k,M)$ then whenever $H$ is a Hilbert space
of dimension $n$, and
$S:H\to X$ and
$T:X\to H$ are bounded operators with $\|T\|\le 1$, we have
$$ |\tr (TS)|\le C\ell(S)(\frac1M+k^{1/2}(\log n)^{-1/2})n^{1/2}.$$
\endproclaim

\demo{Proof}We first choose an orthonormal basis $(e_i)_{i=1}^n$ of
$H$ so that $\|Se_i\|=\|S_{H_i}\|$ where
$H_i=[e_i,e_{i+1},\ldots,e_n].$
Pick $m=[n^{1/2}].$  Then  if $g_1,\ldots,g_n$ are normalized
independent Gaussians,
$$ \ell(S)\ge \bold E(\|\sum_{i=1}^m g_i Te_i\|^2)^{1/2}\ge \|Se_m\|\bold
E(\max_{1\le i\le m}|g_i|).$$
Now (cf. \cite{13} p. 23) this implies that
$$ \|Se_m\| \le C(\log m)^{-1/2}\ell(S)$$ for some universal constant
$C.$  Our choice of $m$ implies that we can replace this estimate by
$$ \|Se_m\| \le C(\log n)^{-1/2}\ell(S)$$ for some universal constant
$C.$

If $E=[e_1,\ldots,e_m]$ then $$|\tr (TS_E)| \le m^{1/2}\ell(S_E) \le
n^{1/4}\ell(S).$$

On the other hand $$|\tr (TS_{H_{m+1}})| \le C \max(M^{-1},k^{1/2} (\log
n)^{-1/2})\ell(S)n^{1/2}$$ by Lemma 2.2.
Combining these results gives us our estimate.\qed\enddemo

\proclaim{Proposition 2.4}Suppose $(X_j)_{j=1}^{\infty}$ is not
sufficiently Euclidean. Then
$\ell_1(X_j)$ is not sufficiently Euclidean.\endproclaim

\demo{Proof}Suppose $(X_j)$ satisfies property $P(k,M).$  Suppose
$n\in\Bbb
N$ and $S:\ell_2^n\to \ell_1(X_j)$ and $T:\ell_1(X_j)\to\ell_2^n$ are any
operators satisfying $TS=I_{\ell_2^n};$ we assume that $\|T\|=1.$ We
write
$S\xi=(S_i\xi)_{i=1}^{\infty}$ and
$T(x_i)_{i=1}^{\infty}=\sum_{i=1}^{\infty}T_ix_i.$

Now $n=\tr(TS)=\sum_{i=1}^{\infty}\tr(T_iS_i).$
On the other hand, by Lemma 2.3 we have that:
$$ |\tr(T_iS_i)|  \le
Cn^{1/2}(\frac 1M +k^{1/2}(\log
n)^{-1/2}\ell(S_i)).\tag 2.1$$

Let $(e_1,\ldots,e_n)$ be any orthonormal basis.  Then by the
Kahane-Khintchine inequality we have that
$$ \ell(S_i) \le C_0\bold E(\|\sum_{j=1}^ng_jS_ie_j\|)\tag 2.2$$
where the $(g_i)_{i=1}^n$ are normalized independent Gaussians, and $C_0$
is a universal constant. Hence
$$ \sum_{i=1}^{\infty}\ell(S_i) \le C_0\bold
E(\|\sum_{j=1}^ng_jSe_j\|)\le C\ell(S)
.\tag 2.2$$

Hence, combining (2.1), (2.2) and (2.3), and since
$n=\sum_{i=1}^{\infty}\tr(T_iS_i),$
$$ n \le C_0n^{1/2}(\frac1M+k^{1/2}(\log n)^{-1/2})\ell(S)\le
C_1n\|S\|(\frac1M+ k^{1/2}
(\log n)^{-1/2})$$
for some universal constant $C_1.$
We thus obtain an estimate
$$ \|S\| \ge  C_1^{-1}\min(M,k^{-1/2}(\log n)^{1/2})$$ for some absolute
constant
$C.$ This establishes the result.
\qed\enddemo

\demo{Remark}In the case when $X=c_0$ we have that $X$ satisfies
$P(ck^{1/2},k)$ for $c>0$ and all $k.$  We thus obtain that $\ell_1(X)$
satisfies $P(c(\log k)^{1/4},k)$ for some $c>0$ and all $k.$  On the
other hand Figiel, Lindenstrauss and Milman \cite{7} established the
upper estimate that $\ell_1(c_0)$ contains a subspace 2-isomorphic to
$\ell_2^k$ which is $(\log k)^{1/2}-$complemented; this estimate is best
possible (see \cite{2}).
 This suggests
our method, while  not optimal, cannot be improved
significantly.\enddemo

 \proclaim{Corollary 2.5}
Suppose $(X_j)_{j=1}^{\infty}$ is not
sufficiently Euclidean. Then
$c_0(X_j)$ is not sufficiently Euclidean.\endproclaim

\demo{Proof}This follows by simple duality.\qed\enddemo

\heading{3. Unconditional bases in $c_0$-products}\endheading

For convenience we define a sequence space $X$ as a K\"othe space  of
real-valued
functions on a countable set $I$ (with counting measure) so that the
canonical basis vectors
$(e_i)_{i\in I}$ form a 1-unconditional basis.  Usually of course we take
$I=\Bbb N$, but for our purposes it is convenienent also to allow $I=\Bbb
N\times \Bbb N$ and certain other alternatives.  A typical element $x$ of
$X$ is of the form $x=(x(i))_{i\in I}.$

Let $(u_n)_{n\in\Cal N}$ be a  set of disjointly supported
vectors in
$X$.  Then $(u_n)_{n\in\Cal N}$ is an unconditional basic sequence, which
is complemented if and only if there exists a biorthogonal sequence
$(u_n^*)_{n\in\Cal N}\in X^*$ with $\supp
u_n^*\subset
\supp u_n$, $u_nu_n^*\ge 0$, $\langle u_n,u_n^*\rangle=1$ and such that
the projection
$$ Px=\sum_{n\in\Cal N}\langle x,u_n^*\rangle u_n$$ is
well-defined and
bounded.  If we define $f_n=u_nu_n^*$ then $f_n\ge 0$, $f_n\in\ell_1(I)$
and $\|f_n\|_1=1$ for all $n\in\Bbb N.$  Under these circumstances we say
that $(u_n)$ is a complemented disjoint sequence and we assume that
$(u_n^*)$ and $(f_n)$ are associated with $(u_n)$.  Note that we can
always replace $u_n$ and $u^*_n$ by $|u_n|$ and $|u^*_n|$ and hence also
assume them positive.

We start with an observation which we will use repeatedly.

\proclaim{Lemma 3.1}Suppose $X$ is a sequence space (modelled on an
index set $I$) and $(u_n)_{n\in\Cal N}$ is a complemented disjoint
sequence.
Let $(A_n)_{n\in\Cal N}$ be any sequence of disjoint sets
such that,
for some $\delta>0$ we have $\|f_n\chi_{A_n}\|_1\ge \delta>0$ for all
$n\in \Cal N.$ Let $v_n=u_n\chi_{A_n}$. Then
$(v_n)_{n\in\Cal N}$ is a complemented
disjoint sequence equivalent to
$(u_n)_{n\in\Cal N}.$ Furthermore the biorthogonal vectors
$v_n^*$ may be chosen so that $v_nv_n^*\le \delta^{-1}f_n.$ \endproclaim

\demo{Proof}Let $P$ be the projection onto $[u_n]_{n\in\Cal N}$ as
defined above.  Let $$v_n^*=\|f_n\chi_{A_n}\|_1^{-1}u_n^*\chi_{A_n}$$ and
define
$$ Qx=\sum_{n\in\Cal N}\langle x,v_n^*\rangle v_n.$$
Then it is easy to verify that $Q$ is a bounded projection,
$Qu_n=v_n$ and
$P(v_n)=\|f_n\chi_{A_n}\|_1u_n.$  This quickly establishes
the equivalence of $(u_n)_{n\in\Cal N}$ and $(v_n)_{n\in\Cal
N}.$  Furthermore $v_n^*v_n\le \delta^{-1}f_n.$
\qed\enddemo

Next suppose $(X_i)_{i=1}^{\infty}$ is a sequence of sequence spaces
modelled on index sets $J_i$ (either finite sets or $\Bbb N)$.  We
suppose that for some $q<\infty$ the spaces $(X_i)$ satisfy a
lower
$q$-estimate uniformly, i.e. there exists $c>0$ so that if $i\in\Bbb N$
and $x_1\ldots,x_n$ are disjoint in $X_i$ then
$$ \|\sum_{k=1}^nx_k\|_{X_i}\ge c(\sum_{k=1}^n\|x_k\|_{X_i}^q)^{1/q}.$$

Let $Y=c_0(X_j)$ be the sequence space on $J=\{(i,j):j\in J_i,\
i\in\Bbb N\}$ of all $x=(x[i])_{i=1}^{\infty}$, where $x[i]\in X_i$ are
so that
$\lim_{i\to \infty}\|x[i]\|_{X_i}=0.$  We define
$$ \|x\|_Y=\max_{i\in \Bbb N}\|x[i]\|_{X_i}.$$

Now suppose $(u_n)_{n\in\Cal N}$ is a complemented disjoint
sequence in
$Y$, with biorthogonal sequence $(u_n^*)$.  As above let
$f_n=u_n^*u_n$. Then define $F_n\in\ell_1(\Bbb N)$ by
$$ F_n(i)=\sum_{j\in J_i}f_n(i,j).$$
We will say that $(u_n)$ is $C$-tempered if
$$ \sum_{i=1}^{\infty}\sup_{n\in\Cal N}F_n(i)\le C.$$

\proclaim{Theorem 3.2}Suppose $(X_i)_{i=1}^{\infty}$ is a sequence of
sequence spaces satisfying a uniform lower $q$-estimate for some
$q<\infty.$  Suppose $(u_n)_{n\in\Cal N}$ is a normalized complemented
disjoint
sequence in $c_0(X_i).$  Then there is a complemented disjoint sequence
$(v_n)_{n\in\Cal N}$ equivalent to $(u_n)_{n\in\Cal N}$ and a partition
$\Cal N=\cup_{n\in\Cal A}\Cal B_n$ of $\Cal N$ with the following
properties:
\newline
(1) For each $i\in\Bbb N$ we have either $\|v_k[i]\|_{X_i}=0$ or
$1$.\newline
(2)
 For some $C$ each $(v_k)_{k\in \Cal
B_n}$ is
$C$-tempered.\newline
(3) There exists an integer $N$ and subsets $(S_n)_{n\in\Cal A}$ of
$\Bbb N$ such that $\sum_{n\in\Cal A}\chi_{S_k}\le N-1$ and $v_k[i]=0$
 whenever $k\in \Cal B_n$ and $i\notin S_n.$
 Hence
 for any finitely nonzero sequence
$(a_n)_{n\in\Cal N}$,
$$ \max_{n\in\Cal A}\|\sum_{k\in\Cal B_n}a_kv_k\|_Y
\le \|\sum_{k\in\Cal N}a_kv_k\|_Y \le
 (N-1)\max_{n\in\Cal A}\|\sum_{k\in \Cal B_n}a_kv_k\|_Y.$$

 \endproclaim

\demo{Proof}As usual we let $P$ be the induced projection on
$[u_n]_{n\in\Cal N}$.

First let $A_n=\{(i,j):\|u_n[i]\|_{X_i}\ge (2\|P\|)^{-1}\}.$  Notice that
$$ \sum_{i\notin A_n}\langle
u_n[i],u_n[i]^*\rangle
\le \frac1{2\|P\|}\sum_{i=1}^{\infty}\|u_n[i]\|_{X_i^*}\le
\frac1{2\|P\|}\|u_n\|_{Y^*}\le \frac12.$$
Hence
$(u_n\chi_{A_n})_{n\in\Cal N}$ is a complemented disjoint sequence
equivalent to $(u_n)_{n\in\Cal N}.$  It follows after some appropriate
renormalization that we can replace $(u_n)_{n\in\Cal N}$ by an equivalent
sequence with the additional property that
$\|u_n[i]\|_{X_i}=1$ or $u_n[i]=0$ for every $i\in\Bbb N.$ For each
$n\in\Cal N$ let $S_n=\{i:\|u_n[i]\|_{X_i}=1\}.$

Next fix any $N\in\Bbb N$ so that $N>1+c^{-q}(1+\|P\|)^q,$ where $c$ is
the constant of the uniform lower $q$-estimate. Let
$\delta=N^{-1}.$  We pick a maximal subset $\Cal A$ of $\Cal N$  with the
property that if $\Cal F$ is a subset of $\Cal A$ with $|\Cal F|\le N$
then
$$ \sum_{i=1}^{\infty}\max_{n\in\Cal F}F_n(i)\ge (1-\delta)|\Cal F|.$$
(Here $F_n(i)=\sum_{j\in J_i}f_n(i,j)=\langle u_n[i],u_n^*[i]\rangle$ as
usual.)

Now let $\Cal F$ be any subset of $\Cal A$ with $|\Cal F|=N.$  We can
partition $\Bbb N$ into $N$ disjoint sets $(A_n)_{n=1}^N$ so that if
$i\in A_n$ then $F_n(i)=\max_{m\in\Cal F}F_m(i).$  Let
$v=\sum_{n\in\Cal F}u_n\chi_{A_n'}$ where $A_n'=\{(i,j):i\in A_n,\ j\in
J_i\}.$  Clearly $\|v\|_Y\le 1.$  However
$$ Pv =\sum_{n\in \Cal F}\sum_{i\in A_n}F_n(i)u_n.$$
Since
$$ \sum_{n\in\Cal F}\sum_{i\in A_n}F_n(i)\ge N(1-\delta)=N-1$$ we
conclude that
$$ \|\sum_{n\in\Cal F}u_n\|_Y \le \|P\|+1.$$
On the other hand,
$$ \|\sum_{n\in\Cal F}u_n\|_Y \ge c\max_{i\in \Bbb N}(\sum_{n\in\Cal F}
\chi_{S_n}(i))^{1/q}.$$
Hence
$$\max_{i\in\Bbb N}\sum_{n\in \Cal F}\chi_{S_n}(i) \le c^{-q}(\|P\|+1)^q
<N.\tag 3.1$$

Now suppose for some $i$ we have $\sum_{n\in\Cal A}\chi_{S_n}(i)\ge N.$
Then we can find a subset $\Cal F$ of $\Cal A$ with $|\Cal F|=N$ so that
$\chi_{S_n}(i)=1$ for $n\in\Cal F$ contradicting (3.1).
We therefore conclude that:
$$ \max_{i\in\Bbb N}\sum_{n\in\Cal A}\chi_{S_n}(i) \le N-1.$$

Now suppose $k\in\Cal N\setminus\Cal A$.  There exists a subset $\Cal F$
of $\Cal A$ with $1\le |\Cal F|\le N-1$ and such that
$$ \sum_{i=1}^{\infty}\max(F_k(i),\max_{n\in\Cal F}F_n(i))<(|\Cal
F|+1)(1-\delta).$$
Hence
$$\align \sum_{i=1}^{\infty}\min(F_k(i),\max_{n\in\Cal
F}F_n(i))&= \sum_{i=1}^{\infty}(F_k(i)+\max_{n\in\Cal
F}F_n(i))-\sum_{i=1}^{\infty}\max(F_k(i),\max_{n\in\Cal F}F_n(i))\\
&> \sum_{i=1}^{\infty}F_k(i)-1+\delta\\
&>\delta.\endalign $$

Thus there exists $n\in\Cal F$ so that
$$ \sum_{i=1}^{\infty}\min(F_k(i),F_n(i))>\delta/N=\delta^2.$$
Now let
$T_k=\{i:F_k(i)<2N^2F_n(i)\}.$  Then $\sum_{i\in
T_k}F_k(i)>\frac12\delta^2.$  In the case when $k\in \Cal A$ we will put
$T_k=S_k.$

We now can partition $\Cal N$ into disjoint sets $(\Cal B_n)_{n\in\Cal
A}$ so that $n\in\Cal B_n$ and if $k\in \Cal B_n$ then $\sum_{i\in
T_k}F_k(i)>\frac12\delta^2$ and $F_k(i)<2N^2F_n(i)$ for $i\in
T_k\subset S_n.$ If we then let $T_k'=\{(i,j):i\in T_k,\ j\in J_i\}$ and
$v_k=u_k\chi_{T_k'}$ then by Lemma
3.1 we have that $(v_k)_{k\in \Cal N}$ is a complemented disjoint
sequence
in $Y$ equivalent to $(u_k)_{k\in\Cal N}.$  Furthermore if $(v_k^*)$ is
the biorthogonal sequence then we have an estimate $v_kv_k^*\le Mf_k$ for
a suitable constant $M.$  It follows that for each $n\in\Cal A$ we have
$$ \langle v_k[i],v_k^*[i]\rangle \le 2MN^2F_n(i)$$
whenever $k\in\Cal B_n$.  Thus the sets $(v_k)_{k\in\Cal B_n}$ are
each $C$-tempered where $C$ is a constant depending only on $c,q$ and
$\|P\|.$

Finally suppose $(a_n)_{n\in\Cal N}$ is finitely non-zero.
Then
$$
\align
 \|\sum_{k\in\Cal N}a_kv_k\|_Y &=\max_{i\in\Bbb N}\|\sum_{k\in\Cal
N}a_kv_k[i]\|_{X_i}\\
&=\max_{i\in\Bbb N} \|\sum_{n\in\Cal A,i\in S_n}\sum_{k\in
\Cal B_n}a_kv_k[i]\|_{X_i}\\
&\le (N-1)\max_{n\in\Cal A}\|\sum_{k\in\Cal B_n}a_kv_k\|_Y.
\endalign
$$
This completes the proof.\qed\enddemo

Let us first use this result to give a simpler proof of the result of
\cite{3} that $c_0(\ell_1)$ has a unique unconditional basis (up to
permutation).

\proclaim{Theorem 3.3}The space $c_0(\ell_1)$ has a unique unconditional
basis.\endproclaim

\demo{Proof}We start with the remark that $c_0(\ell_1)$ is not
sufficiently Euclidean (cf. Bourgain \cite{2} or Corollary 2.5 above).
Hence any complemented unconditional basic sequence is equivalent to a
complemented positive disjoint sequence in $c_0(\ell_1)^m$ for some $m$.
We can clearly suppose $m=1.$

We next show that any $C$-tempered
$C$-complemented
disjoint sequence $(u_n)_{n\in\Cal N}$ is $K$-equivalent to the
standard
$\ell_1-$basis where
$K$ depends only on $C.$ Indeed we may suppose $\|u_n[i]\|_{1}=1$ or
$u_n[i]=0$ for each $i,n$.  Let $G(i)=\max_n\langle
u_n[i],u_n^*[i]\rangle.$ Then $\|G\|_1 \le C$.

Now
$$\align
 \|\sum_{n\in\Cal N}a_nu_n\|_{c_0(\ell_1)}&\ge \frac1C
\sum_{i=1}^{\infty}G(i)\sum_{u_n[i]\neq 0}|a_n| \\
&= \frac1C\sum_{n\in\Cal N}|a_n|\sum_{u_n[i]\neq 0}G(i) \\
&\ge \frac1C\sum_{n\in\Cal N}|a_n|\sum_{i=1}^{\infty}\langle
u_n[i],u_n^*[i]\rangle\\
&=\frac1C\sum_{n\in\Cal N}|a_n|.
\endalign
$$

It follows that any unconditional basis of $c_0(\ell_1)$ is equivalent to
the canonical unconditional basis of $c_0(X_n)$ where each $X_n$ is
either $\ell_1$ or $\ell_1^m$ for some $m=m(n).$  However there must be
infinitely many indices $n$ for which $X_n=\ell_1$ (since $c_0(\ell_1)$
cannot be decomposed as $\ell_1\oplus Z$ where $Z$ contains no copy of
$\ell_1.$)  It then easily follows that $c_0(\ell_1)$ has a unique
unconditional basis.\qed\enddemo

\proclaim{Theorem 3.4}Suppose $1<p_n<\infty$ and $p_n\downarrow 1.$
Let
$(N_n)$ be an increasing sequence of natural numbers such that
$(p_n-p_{n+1})\log N_n$ is bounded and $N_{n+1}\ge 2N_n.$  Then
$c_0(\ell_{p_n}^{N_n})$ has a unique unconditional basis.\endproclaim

\demo{Proof}We suppose the sequence $(p_n)$ fixed and first consider
$c_0(\ell_{p_n}^{M_n})$ for any sequence of integers $(M_n).$  It is easy
to see by considering the ultraproduct $\prod_{\Cal U}\ell_{p_n}^{M_n}$
(which is an $L_1-$space) that the sequence $(\ell_{p_n}^{M_n})$ is not
sufficiently Euclidean.  By Corollary 2.5 the space
$Y=c_0(\ell_{p_n}^{M_n})$
is also not sufficiently Euclidean and therefore  every complemented
unconditional basic sequence in $Y$ is equivalent to a complemented
disjoint sequence in $Y^r=c_0(\ell_{p_n}^{rM_n})$ for some $r\in \Bbb
N.$

In this case our spaces $X_n$ are modelled on the sets
$J_n=\{1,2,\ldots,M_n\}.$  Thus $c_0(\ell_{p_n}^{M_n})$ is modelled on
the set $\{(n,k):\ n\in\Bbb N,\ 1\le k\le M_n\}.$

Now suppose $(u_n)_{n\in\Cal N}$ is a $C$-tempered $C$-complemented
disjoint sequence in $c_0(\ell_{p_n}^{M_n})$ with the property that
$\|u_n[i]\|_{\ell_{p_i}^{M_i}}=1$ or $0$ for each $i$, and let
$A=\{i:\exists n,\ u_n[i] \neq 0\}.$
We claim that:

{\it Claim: there exists a constant $K=K(C),$ an integer $r$ depending
only on
$C$,
a subset $B$ of $A$ with $|B|\le r$ and $P_i\le r M_i$ for $i\in B$ so
that
$(u_n)_{n\in\Cal N}$ is $K$-equivalent to the canonical basis of
$(\sum_{i\in B}\oplus\ell_{p_i}^{P_i})_{c_0}.$}

To this end let $G(i)=\sup \langle u_n[i],u_n^*[i]\rangle$ as usual.  We
have $\sum_{i=1}^{\infty}G(i)\le C.$
We first consider the case when $G(i)\le \frac1{20}$ for every $i$.
Then it is possible to find a finite increasing sequence of integers
$(k_j)_{j=0}^N$ where $N=[10C]$ depends on $C$ such that
$$ \sum_{i<k_j}G(i) \le j/10$$ and
$$ \sum_{i\le k_j}G(i)\ge j/10.$$
Notice that
$$ \sum_{i>k_N}G(i)\le 1/10.$$
It follows that for each $n$ there exist at least three values of $1\le
j\le N$ so that
$$ \sum_{i=k_{j-1}+1}^{k_j}\langle u_n[i],u_n^*[i]\rangle>\frac{1}{4N}.$$
We can then assign to each $n$ a value of $j$ which is neither the
largest or smallest with this property.  In this way we partition $\Cal
N$ into sets $(\Cal N_j)_{2\le j\le N-1}.$

Consider $(u_n)_{n\in\Cal N_j}.$  We note that this is equivalent (with
constants depending only on $C$) to each of $(v_n)_{n\in\Cal N_j}$ and
$(w_n)_{n\in\Cal N_j}$ where $v_n[i]=u_n[i]$ if $i\le k_{j-1}$ and $0$
otherwise while $w_n[i]=u_n[i]$ if $i>k_j$ and $0$ otherwise.

Now for any finitely non-zero sequence $(a_n)_{n\in\Cal N_j}$ we have
$$ \|\sum_{n\in\Cal N_j}a_nv_n\|_Y\le (\sum_{n\in\Cal
N_j}|a_n|^{q_{j-1}})^{1/q_{j-1}}$$ where $q_j=p_{k_j}.$  On the other
hand
$$ \|\sum_{n\in\Cal N_j}a_nw_n\|_Y\ge \max_{i>k_j}(\sum_{u_n[i]\neq
0}|a_n|^{q_j})^{1/q_j}.$$
 Thus
$$ \|\sum_{n\in\Cal N_j}a_nw_n\|_Y^{q_j} \ge \frac1C
\sum_{i>k_j}\sum_{u_n[i]\neq 0}G(i)|a_n|^{q_j}.$$

Now for fixed $n$
$$ \sum_{i>k_j,u_n[i]\neq 0}G(i) \ge \sum_{i>k_j}\langle
u_n[i],u_n^*[i]\rangle \ge \frac1{4N}.$$
It follows that $(u_n)_{n\in\Cal N_j}$ satisfies an upper
$q_{j-1}-$estimate and a lower $q_j-$estimate with constants depending
only on $C.$  We next estimate $|\Cal N_j|.$  In fact
$$
\align
 \frac1{4N}|\Cal N_j|&\le \sum_{n\in\Cal
N_j}\sum_{i=k_{j-1}+1}^{k_j}\langle u_n[i],u_n^*[i]\rangle \\
&\le \sum_{i=k_{j-1}+1}^{k_j}G(i)M_i\\ & \le C\max_{k_{j-1}<i\le k_j}M_i
\endalign
$$

Hence if we select $k_{j-1}<i\le k_j$ appropriately we have
$(u_n)_{n\in\Cal N_j}$ equivalent to a subset of the standard basis of
$\ell_{p_i}^{rM_i}$ where $r$ and the constant of equivalence depend only
on $C.$

We must now treat the case when $G(i)>\frac1{20}$ for some $i.$  In this
case we split $\Cal N$ into two groups $\Cal N'$ and $\Cal N''$ where
$\Cal N'=\{n:\exists i, \ \langle u_n[i],u_n^*[i]\rangle >\frac1{20}\}$
and
$\Cal N''$ is the remainder.  Then $\Cal N''$ can be treated as before.
For $\Cal N'$ we note that $(u_n)_{n\in\Cal N'}$ is equivalent to a
sequence $(u'_n)$ where $u'_n[i]=u_n[i]$ for precisely one index $i=i_n$
such that $\langle u_n[i_n],u_n^*[i_n]\rangle>1/20$ and is zero
elsewhere. The appropriate representation of
$(u_n)_{n\in\Cal N'}$ follows once we observe that the set
$\{i_n:n\in\Cal N'\}$ is bounded in cardinality with a bound depending
only on $C.$ But this is clear since $G(i_n)>1/20$ but $\sum G(i)\le C.$

Thus the claim is established.

Returning to our original hypotheses we see that if $(u_n)$ is any
unconditional basis of $c_0(\ell_{p_n}^{N_n})$ then $(u_n)$ is equivalent
to the canonical basis of $c_0(\ell_{p_n}^{M_n})$ where $M_n\le rN_n$ for
all $n$ and some fixed $r.$  By the same token the canonical basis of
$c_0(\ell_{p_n}^{N_n})$ is equivalent a subset of $(u_n)^s$ for some
$s\in\Bbb N.$

Now the additional hypotheses on $N_n$ ensure that the original basis is
equivalent to its square.  Hence the $s$-fold product $(u_n)^s$ is
equivalent to a subset
of the canonical basis and so it follows from the Cantor-Bernstein
principle (apparently first noticed by Mityagin, \cite{14},\cite{16} and
\cite{17}), that
$(u_n)^s$ and the original basis are equivalent.

Thus the canonical bases of $c_0(\ell_{p_n}^{sM_n})$ and
$c_0(\ell_{p_n}^{sN_n})$ are equivalent.
Let $\Cal M=\{(i,j):1\le j\le M_i\}$ and $\Cal N=\{(i,j):1\le j\le
N_i\}.$
Let us suppose the former basis is indexed by $\Cal M^s=\{(i,j): 1\le
j\le sM_i\}$ and the latter by $\Cal N^s=\{(i,j): 1\le j\le sN_i\}.$
Let $\varphi:\Cal M^s\to\Cal N^s$ be a bijection implementing the claimed
equivalence of bases.  By elementary considerations concerning $c_0-$sums
it is clear that for each fixed $i$ the set $\{\varphi(i,j): 1\le j\le
M_i\}$ can have
at most $t$ distinct first co-ordinates where $t$ depends only on the
constant of equivalence; similarly for each fixed $i$ the set of
possible first co-ordinates of $\varphi^{-1}(i,j)$ can be bounded by the
same $t.$  For each $(a,b)\in \Bbb N^2,$ let $E_{ab}$ be the set of
$(i,j)$ so that $i=a$ and the first co-ordinate of $\varphi(i,j)$ is $b$.
Then let $F_{ab}$ be a subset of $E_{ab}$ of size $[|E_{ab}|/s].$
Restricting $\varphi$ to $\cup_{(a,b)}F_{ab}$ produces
 an equivalence between the bases of
$ c_0(\ell_{p_n}^{M_n-\alpha_n})$ and $c_0(\ell_{p_n}^{N_n-\beta_n})$
where
$0\le \alpha_n,\beta_n\le (s-1)t$ for each $n.$  This clearly implies
the equivalence of
$(u_n)$ and the original basis.\qed\enddemo

\heading{4. Uniqueness of unconditional bases in
$c_0$-products of right-dominant spaces}\endheading

We first introduce some standard notation.  Let $A,B$ be subsets of
$\Bbb N.$ We write $A<B$ to indicate that $\max\{a:a\in A\}<\min\{b:b\in
B\}.$

Let $X$ be a sequence space modelled on $\Bbb N.$
 We say that $X$ is
{\it right-dominant} if there is a constant $\kappa=\kappa(X)$ so that
whenever $u_1,\ldots,u_n$ and $v_1,\ldots,v_n$ are any two disjointly
supported sequences satisfying $\supp u_k<\supp v_k$ and
$\|u_k\|_X=\|v_k\|_X$ for $1\le k\le n$ then $\|\sum_{k=1}^nu_k\|_X\le
\kappa \|\sum_{k=1}^nv_k\|_X.$
 We say that $X$ is
{\it left-dominant} if there is a constant $\rho=\rho(X)$ so
that
whenever $u_1,\ldots,u_n$ and $v_1,\ldots,v_n$ are any two disjointly
supported sequences satisfying $\supp u_k<\supp v_k$ and
$\|u_k\|_X=\|v_k\|_X$ for $1\le k\le n$ then $\|\sum_{k=1}^nv_k\|_X\le
\rho \|\sum_{k=1}^nu_k\|_X.$

Left and right-dominant spaces were studied in \cite{4}.  It is
established (Lemma 5.2 of \cite{4}) that in these spaces there is exactly
one
$r=r(X)$
(the {\it index} of $X$) so that $\ell_r$ is disjointly finitely
representable in $X$.  If $X$ is right-dominant then $X$ satisfies an
upper $r$-estimate and a lower $s$-estimate for any $s>r$; the
corresponding dual statements hold for left-dominant spaces.  Clearly if
a space $X$ is both left and right-dominant then $X=\ell_r.$

\proclaim{Theorem 4.1}Let $X$ be a right-dominant sequence space with
$r(X)=1.$  Then every complemented unconditional basic sequence in
$c_0(X)$ is equivalent to a subsequence of the canonical
basis.\endproclaim

\demo{Remark}In particular this applies when $X$ is a Nakano space
$\ell_{(p_n)}$ where $p_n\downarrow 1$ or when $X$ is Tsirelson space $T$
(see \cite{4}).\enddemo

\demo{Proof}In this case we note that in the notation of Section 3,
$J_i=\Bbb N$ for all
$i\in\Bbb N.$ We first note that by Corollary 2.5,
$c_0(X)$ is not sufficiently Euclidean. Hence by \cite{4} Theorem 3.5,
every complemented
unconditional basic sequence is equivalent to a complemented positive
disjoint sequence in $c_0(X)^N$ for some $N$ and hence also to a
complemented disjoint sequence in $c_0(X).$

Now by Theorem 3.2  it will suffice to show that if $(u_n)_{n\in\Cal N}$
is a
$C$-tempered $C$-complemented unconditional basic sequence then
$(u_n)_{n\in\Cal N}$
is $K$-equivalent to a subsequence of the canonical basis of $c_0(X)$
where
$K$ depends only on $C.$  In fact we will show that it is $K$-equivalent
to a subsequence of the canonical basis of $\ell_{\infty}^N(X)$ where $N$
depends only on $C.$

We may suppose that, as before, $\|u_n[i]\|_X=1$ or $u_n[i]=0.$  Let
$f_n=u_nu_n^*$ and $F_n(i)=\langle u_n[i],u^*_n[i]\rangle.$  Let
$G(i)=\max_n F_n(i)$ so that $\|G\|_1\le C.$

Pick any integer $N>2C.$  For each $n\in\Cal N$ we pick natural
numbers
$r_n\le s_n$ such that
$$ \sum_{i=1}^{\infty}\sum_{j=1}^{r_n-1}f_n(i,j) <\frac14$$
$$ \sum_{i=1}^{\infty}\sum_{j=1}^{r_n}f_n(i,j) \ge\frac14$$
$$ \sum_{i=1}^{\infty}\sum_{j=1}^{s_n-1}f_n(i,j) <\frac34$$
$$ \sum_{i=1}^{\infty}\sum_{j=1}^{s_n}f_n(i,j) \ge\frac34.$$
We will argue by Hall's Marriage Lemma (see Bollobas \cite{1}) that it
is possible to find an
map $\varphi:\Cal N\to\Bbb N$ such that $\varphi(n)\in
[r_n,s_n],$ with $|\varphi^{-1}(k)|\le N$ for all $k\in\Bbb N.$
Indeed if not, the Marriage Lemma implies there is a minimal finite
subset
$\Cal M$ of
$\Cal N$ such that $N|\cup_{n\in\Cal M}[r_n,s_n]| < |\Cal M|.$  It
follows
easily from the mininality that $\cup_{n\in\Cal M}[r_n,s_n]$ is an
interval $[a,b].$ From the disjointness of the $(f_n)$ we have that
$$ \sum_{j=a}^b\sup_{n\in\Cal M}f_n(i,j)\le (b-a+1)G(i)$$
so that
$$ \sum_{i=1}^{\infty}\sum_{j=a}^bf_n(i,j) \le (b-a+1)C$$
However
$$ \sum_{i=1}^{\infty}\sum_{j=a}^bf_n(i,j) \ge \frac12|\Cal M|$$
so that $|\Cal M| \le 2(b-a+1)C<N(b-a+1)$ which is a contradiction.

We can now split $\Cal N$ into at most  $N$ disjoint subsets
$(\Cal N_k)_{k
\in\Cal M}$ so that $\varphi$ is injective on each $\Cal N_k.$

For each $n\in\Cal N$ let $A_n=\{(i,j):j\le\varphi(n)\}$ and
$B_n=\{(i,j):j\ge\varphi(n)\}.$  Then, by Lemma 3.1, we have that
$(u_n)_{n\in\Cal N}$ is equivalent to both $(u_n\chi_{A_n})_{n\in\Cal N}$
and
$(u_n\chi_{B_n})_{n\in\Cal N}$ with constants of equivalence depending
only on
$C.$

Now suppose $k\in\Cal M$, and let $(a_n)_{n\in\Cal N_k}$ be a finitely
non-zero sequence.
Then, by the right-dominance property, for each $i$
$$ \|\sum_{n\in\Cal N_k}a_nu_n\chi_{A_n}[i]\|_X \le
\kappa\|\sum_{n\in\Cal N_k}a_n\|u_n\chi_{A_n}[i]\|_Xe_{\varphi(n)}\|_X$$
Hence
$$ \|\sum_{n\in\Cal N_k}a_nu_n\chi_{A_n}[i]\|_Y \le
\kappa\|\sum_{n\in\Cal N}a_ne_{\varphi(n)}\|_X.\tag 4.1$$

In the opposite direction, again by the right-dominance property, we have
$$ \|\sum_{n\in\Cal N_k}a_n\|u_n\chi_{B_n}[i]\|_X e_{\varphi(n)}\|_X \le
\kappa \|\sum_{n\in\Cal N_k}a_nu_n\chi_{B_n}[i]\|_X.\tag 4.2$$

Combining (4.1) and (4.2),
$$ \|\sum_{n\in\Cal N_k}a_n\sum_{i=1}^{\infty}G(i)\|u_n\chi_{B_n}[i]
\|_Xe_{\varphi(n)}\|_X \le C\kappa\|\sum_{n\in\Cal N_k}a_nu_n\|_Y.\tag
4.3$$

For each $n$ let $I_n=\{i: 8CG(i)\ge \|u_n^*[i]\|_{X^*}\}.$
Then
$$ \sum_{i\notin I_n}G(i)\le
\frac{1}{8C}\sum_{i=1}^{\infty}\|u_n^*[i]\|_{X^*} \le \frac18.$$
Hence
$$ \sum_{i\notin I_n}\langle u_n[i]\chi_{B_n},u_n^*[i]\rangle \le
\frac18.$$
However by choice of $\varphi(n)$ we have
$$ \sum_{i=1}^{\infty}\langle u_n[i]\chi_{B_n},u_n^*[i]\rangle \ge
\frac14.$$
Thus we have
$$
\align
\frac18 &\le \sum_{i\in I_n}\|u_n[i]\chi_{B_n}\|_X\|u_n^*[i]\|_{X^*}
\\
&\le 8C\sum_{i\in I_n} G(i)\|u_n[i]\chi_{B_n}\|_X\\
&\le 8C\sum_{i=1}^{\infty} G(i)\|u_n[i]\chi_{B_n}\|_X.\endalign$$
The estimate above combined with (4.3) yields the inequality
$$ \|\sum_{n\in\Cal N_k}a_ne_{\varphi(n)}\|_X \le 64C^2\kappa
\|\sum_{n\in\Cal N_k}a_nu_n\|_Y.$$

Thus each $(u_n)_{n\in\Cal N_k}$ is equivalent to $(e_{\varphi(n)})_{n\in
\Cal N_k}$ in $X$ with constant of equivalence depending only on $C$.
Since $|\Cal M| \le N$ where $N$ depends only on $C$, the result is
proved. \qed\enddemo

Let us say that an unconditional basis $(u_n)_{n\in\Cal N}$ is {\it
molecular} if there exists a constant $C$ and a natural
number
$N$ so
that if $\Cal N$ is partitioned into $N$ disjoint sets $(\Cal
N_k)_{k=1}^N$ then there exists a proper subset $\Cal M$ of
$\{1,2,\ldots,N\}$ such that $(u_n)_{n\in\Cal N}$ is $C$-equivalent
to a subset of $\cup_{k\in\Cal M}(u_n)_{n\in\Cal N_k}.$
Otherwise we will say that $(u_n)_{n\in\Cal N}$ is {\it
non-molecular}. It follows from the quantitative form of the
Cantor-Bernstein
principle \cite{14,16,17} that $(u_n)_{n\in\Cal N}$ is molecular if and
only if there is  a constant $C$ so that if $\Cal N$ is partitioned into
$N$ disjoint sets $(\Cal N_k)_{k=1}^N$ then there is a proper subset
$\Cal M$ of $\{1,2,\ldots,N\}$ so that $(u_n)_{n\in \Cal N_k,k\in\Cal
M}$ is $C$-equivalent to $(u_n)_{n\in\Cal N}.$  Let us note that any
subsymmetric basis is molecular with $N=2$ as is the usual basis of
$(\sum_{n=1}^{\infty}\ell_q^n)_{\ell_p}$ when $1\le p,q<\infty.$  The
canonical basis of $\ell_p\oplus\ell_q$ for $p\neq q$ is molecular with
$N=3.$

\proclaim{Lemma 4.2}Suppose $(u_n)_{n=1}^{\infty}$ is a
non-molecular unconditional basis.  Then for any
$\epsilon>0$,
$N\in\Bbb N$ and constant $C$ there exists $M>N$ and subsets
$(A_k)_{k=1}^{M}$ of $\Bbb N$ so that:\newline
(1) For any subset $\Cal M$ of $\{1,2,\ldots, M\}$ with $|\Cal M|<N$
then $(u_n)_{n=1}^{\infty}$ is not $C$-equivalent to any subset of
$\{u_n: n\in \cup_{k\in\Cal M}A_k\},$ and\newline
(2) $\frac1M\sum_{k=1}^M\chi_{A_k}\ge (1-\epsilon)\chi_{\Bbb N}.$
\endproclaim

\demo{Proof}It suffices to consider the case when $\epsilon=r/s$ is
rational.  We then may pick an integer $m$ so large that $mr>N$ and
so that if $L=\binom{ms}{mr}$ then we can partition $\Bbb N$ into $L$
sets so that $(u_n)$ is not $C$-equivalent to a subset of
$(u_n)_{n\in \Cal N}$ where $\Cal N$ is the union of any $L-1$ sets.

Let $\Omega$ be the collection of all $m(s-r)$ subsets of
$\{1,2,\ldots,ms\}.$  We can partition $\Bbb
N=\cup_{\omega\in\Omega}B_{\omega}$ so that $(u_n)_{n\in\Bbb N}$ is not
$C$-equivalent to a subset of $(u_n)_{n\in\Cal N}$ where $\Cal
N=\cup_{\omega\in D}B_{\omega}$ for some proper subset $D$ of $\Omega.$

Now let $A_k=\cup_{\omega:k\in\omega}B_{\omega}$ for $1\le k\le M=ms.$ It
is clear that
$$ \sum_{k=1}^{ms}\chi_{A_k}=m(s-r)\chi_{\Bbb N}$$ so that (2) holds.
Suppose $(u_n)_{n\in\Bbb N}$ is $C$-equivalent to a subset of
$(u_n)_{n\in \Cal N}$ where $\Cal N=\cup_{k\in\Cal M}A_k.$  Then we
have $\cup_{k\in\Cal M}A_k=\Omega$ whence
$|\Cal M| >mr.$  \qed\enddemo

\proclaim{Theorem 4.3}Let $X$ be a space with nontrivial cotype and
an unconditional basis
$(u_n).$ If
$c_0(X)$ has a unique unconditional basis then $(u_n)$ is
molecular.
\endproclaim

\demo{Proof} We will assume, on the contrary that the basis $(u_n)$ is
not molecular. Let us regard
$X$ as a sequence space so that
the given unconditional basis is identified with $(e_n)_{n\in\Bbb N}.$
We start by using Lemma 4.2 repeatedly to generate for each $r\in\Bbb
N$, subsets $(A_{rk})_{k=1}^{M_r}$ of $\Bbb N$ so that:\newline
(1) For any subset $\Cal M$ of $\{1,2,\ldots,M_r\}$ with $|\Cal M|<r$ the
basis $(e_n)_{n=1}^{\infty}$ is not $r$-equivalent to any subset of
$\{e_n:n\in\cup_{k\in\Cal M}A_{nk}\}$, and \newline
(2) $M_r^{-1}\sum_{k=1}^{M_r}\chi_{A_{rk}}\ge (1-2^{-(r+1)})\chi_{\Bbb
N}.$

Now for each $s\in\Bbb N$ let $P_s=\prod_{r=1}^sM_r$ and let
$(B_{sk})_{k=1}^{P_s}$ be a listing of all sets of the form
$\cap_{r=1}^s A_{rk_r}.$  We observe that
$$ P_s^{-1}\sum_{k=1}^{P_s}\chi_{B_{sk}}\ge \frac12\chi_{\Bbb
N}.$$

Consider the index set $I=\{(s,k):\ 1\le k\le P_s,\ s\in\Bbb N\}.$  We
will consider the space $c_0(X)$ as a sequence space modelled on $I\times
\Bbb N$.

Consider now the block basic sequence
$$ u_{sn}=\sum_{k=1}^{P_s}e_{skn}.$$
If we define the biorthogonal functionals
$$ u^*_{sn}=\frac{1}{P_s}\sum_{k=1}^{P_s}e_{skn}$$ then it is clear that
$(u_{sn})_{s,n}$ is a complemented disjoint sequence
equivalent to the canonical basis of $c_0(X).$

Now let $D=\{(s,k,n):\ n\in B_{sk}\}.$  Then
$\|u_nu_n^*\chi_D\|_1\ge \frac12.$  It follows from Lemma 3.1 that
$(u_{sn}\chi_D)_{s,n}$ is also a complemented disjoint  sequence
equivalent to the canonical basis of $c_0(X).$

The basis vectors $(e_{skn})$ for $(s,k,n)\in D$ span a complemented
subspace $Y$ of $c_0(X)$ which by the above remark contains a
complemented copy of $c_0(X).$  By the Pe\l czy\'nski decomposition
argument, $Y$ is isomorphic to $c_0(X).$  If we assume that $c_0(X)$ has
a unique unconditional basis then it will follow that the whole basis
$(e_{skn})_{(s,k)\in I,n\in\Bbb N}$ is $C$-equivalent, for some $C$, to
its subset $(e_{skn})_{
(s,k,n)\in D}.$

Thus we can partition $D$ into subsets $(D_t)_{t=1}^{\infty}$ so that
each subset $(e_{skn})$ for $(s,k,n)\in D_t$ is $C$-equivalent to the
canonical basis $(e_n)$ of $X$ while any subset obtained by picking one
element from each $D_t$ is $C$-equivalent to the standard $c_0$-basis.
>From this and the fact that $X$ has a lower-estimate it is clear that for
fixed $(s,k)$ at most finitely many $D_t$ can intersect the set of all
$(s,k,n)$ for $n\in\Bbb N$.
Note also that the set of $(s,k)$ such that $(s,k,n)\in D_t$ for some
$n$ must also be uniformly bounded by some constant $K$ again by
the lower estimate on
$X$.

 In particular for any $s_0$ there exists $t$
so that if $(s,k,n)\in D_t$ then $s>s_0.$  Hence, the canonical basis of
$X$ is $C$-equivalent to a subset of $\cup_{(s,k)\in\Cal M}B_{sk}$ where
$(s,k)\in\Cal M$ implies $s>s_0$ and $|\Cal M|\le K.$  Now each $B_{sk}$
is contained in some $A_{s_0,k}$ and so we must have $K\ge s_0.$
By choosing $s_0$ large enough we get a contradiction.\qed\enddemo

We now state without proof a  general theorem which can
be proved by exactly the same argument.

\proclaim{Theorem 4.4}Suppose $1\le p<\infty$ and suppose $X$ is a Banach
space with a non-molecular unconditional basis $(u_n)_{n=1}^{\infty}$
with
the property that it does not contain subsets unifomrly equivalent to the
unit vector bases of $\ell_p^m$ for $m=1,2,\ldots.$  Let
$(u_{mn})_{m,n=1}^{\infty}$ be the induced basis of $\ell_p(X).$  Then
there is a subset $\Cal A\subset \Bbb N\times \Bbb N$ so that
$(u_{mn})_{(m,n)\in\Cal A}$ is non-equivalent to the full basis
$(u_{mn})$ and spans a subspace isomorphic to $\ell_p(X).$\endproclaim

We conclude with a theorem which gives us a large number of examples of
right-dominant spaces with non-molecular unconditional bases.

\proclaim{Theorem 4.5}Suppose $X$ is a right-dominant sequence space
with
$r(X)=r.$  Suppose the canonical basis is molecular.  Then
$X=\ell_r.$\endproclaim

\demo{Proof}It is enough to show that $X$ is left-dominant.  Let us
assume the contrary.  Then we claim:
\newline
{\it Claim: Given any $a\in\Bbb N$ and $C>0$ there exists $b>a$ so that
$(e_k)_{a<k\le b}$ is not $C$-equivalent to any subset of $(e_k)_{k\le
a}\cup(e_k)_{b<k}.$}

To prove the claim let $C_1>C^2\kappa+a.$  Since $X$ is not
left-dominant
there exist disjoint sequences $(u_n)_{n=1}^N$ and $(v_n)_{n=1}^N$ with
finite supports so that $a<\supp u_n<\supp v_n$ for each $n$,
$\|u_n\|_X=\|v_n\|_X$ and
$$ \|\sum_{n=1}^Nv_n\|_X>C_1\|\sum_{n=1}^Nu_n\|_X.$$
Pick $b$ so large that $\supp v_n\le b$ for all $n.$  Suppose
$(e_k)_{a<k\le b}$ is $C$-equivalent to some subset of $(e_k)_{k\le
a}\cup (e_k)_{k>b}.$  Then there exist $(w_n)_{n=1}^N$ each with finite
disjoint supports not intersecting $(a,b]$ so that $\|w_n\|_X=\|u_n\|_X$
for
$1\le n\le N$ and
$$ \|\sum_{n=1}^Nw_n\|_X \le C^2\|\sum_{n=1}^Nu_n\|_X.$$
Let $\Cal M=\{n:\ \supp w_n\cap [1,a]\neq \emptyset\}.$  Then $|\Cal
M|\le a.$
Thus
$$ \|\sum_{n\in\Cal M}v_n\|_X\le a\|\sum_{n=1}^Nu_n\|_X.$$
On the other hand
$$ \|\sum_{n\notin\Cal M}v_n\|_X\le \kappa \|\sum_{n\notin\Cal
M}w_n\|_X\le C^2\kappa \|\sum_{n=1}^Nu_n\|_X.$$
It follows that
$C_1< C^2\kappa +a$ contrary to assumption. This establishes the claim.

To prove the theorem we use the claim to find an increasing sequence
$(a_n)_{n=1}^{\infty}$ so that $(e_k)_{a_n< k\le a_{n+1}}$ is not
$n$-equivalent to any subset of $(e_k)_{k\le a_n}\cup(e_k)_{k>a_{n+1}}.$
Then fix any $s\in \Bbb N$ and consider the sets
$A_j=\cup\{(a_n,a_{n+1}]:\ n\cong j \mod s\}$ for $0\le j\le s-1.$
Now the sets $(e_n)_{n\in A_j}$ partition the basis into $s$ sets in such
a way that no
$(s-1)$ sets contain a subset equivalent to the original
basis. This contradicts our assumption that the basis is
molecular.\qed\enddemo

\demo{Examples}We can now give many examples of spaces $X$ with a unique
unconditional basis but such that the $c_0$-product $c_0(X)$ fails to
have unique unconditional basis. This will answer negatively a question
raised in \cite{3}.

 In fact if $X$ is right-dominant and
$c_0(X)$ has unique unconditional basis then $X$ must be one of the
three spaces $c_0,\ell_1$ or $\ell_2.$  This follows by observing that
if it is not in this list then $r(X)<\infty$ and hence $X$ has cotype.
Then
the preceding Theorems 4.3 and 4.5 show that $X=\ell_r$ for some finite
$r$.  The uniqueness then forces either $r=1$ or $r=2.$

On the other hand there are many known examples of right-dominant spaces
with unique unconditional bases.  In \cite{3} 2-convexified
Tsirelson space is shown to have unique unconditional bases.  In
\cite{4} Tsirelson space itself and certain Nakano spaces $\ell_{(p_n)}$
are shown to have unique unconditional bases.  These latter examples
satisfy $r(X)=1$ so that we can apply Theorem 4.1.  The second
non-equivalent basis constructed in Theorem 4.3 is indeed equivalent to
a subset of the original basis.\qed\enddemo

\enddocument
\bye